\documentclass[conference, letterpaper, 10pt]{ieeeconf}
\IEEEoverridecommandlockouts
\usepackage{amsfonts,amsmath,amssymb}
\usepackage{times}
\usepackage{bm}
\usepackage{graphicx}
\usepackage{balance}
\usepackage{setspace}
\usepackage{balance}
\usepackage{comment}
\usepackage{cases}
\graphicspath{{./fig_eps/}}
\usepackage{here}
\usepackage{ascmac} 
\usepackage{array}
\usepackage{float}
\usepackage[linesnumbered, ruled]{algorithm2e}
\usepackage{cite}
\usepackage{url}
\usepackage{subfigure}

\newtheorem{mythm}{Theorem}
\newtheorem{mypro}{Problem}

\newtheorem{myas}{Assumption}

\newcommand{\rfig}[1]{Fig.\,\ref{#1}} 
 
\newcommand{\req}[1]{\eqref{#1}}

\newcommand{\rsec}[1]{Section\,\ref{#1}}
\newcommand{\rpro}[1]{Problem\,\ref{#1}}

\newcommand{\ras}[1]{Assumption\,\ref{#1}}

\newcommand{\ralg}[1]{Algorithm\,\ref{#1}}

\newcommand{\qedwhite}{\hfill \ensuremath{\Box}}
\begin{document}
\title{\LARGE \textbf{Safe and active parameter exploration for event-triggered control}} 
\author{Kazumune~Hashimoto 
\thanks{Kazumune Hashimoto is with the Graduate School of Engineering Science, Osaka University, Osaka, Japan (e-mail: kazumune.hashimoto@hopf.sys.es.osaka-u.ac.jp). 
This work is supported by ERATO HASUO Metamathematics for Systems Design Project (No. JPMJER1603), JST.}}
\maketitle

\begin{abstract}
This paper presents a framework of learning parameter space for event-triggered control. In particular, our goal is to find a set of parameters for the event-triggered condition, such that certain specifications on safety and convergence properties are satisfied. The exploration strategy is based on the Gaussian process-based active learning, in which, for each iteration, the parameter with having the largest variance is evaluated. Moreover, we provide a theoretical analysis, so that the derived parameter space satisfies both convergence and safety. Finally, a numerical simulation is given to illustrate the effectiveness of the approach. 
\end{abstract}

\begin{keywords}
Event-triggered control, Networked control systems, Gaussian process regression. 
\end{keywords}

\section{Introduction}
Resource-aware control, such as event-triggered control and self-triggered control, has attracted much attention in recent years as one of the ways to reduce the \textit{network resource utilization} in networked control systems (NCSs), such as communication bandwidth and energy consumption of the battery powered devices. In the event and self-triggered control, transmissions of sensor measurements and the control updates are not given periodically, but are given only when they are necessary. For example, the event-triggered control executes the control updates when some events based on the current state measurements are triggered. 
Until now, event and self-triggered control have been investigated for various types of formulations. 
Early works aim at designing the event/self-triggered control based on analyzing Lyapunov-based stability, see, e.g., \cite{tabuada2007a,heemels2011a,dolk2017c,mazo2010}. 
More recently, periodic event-triggered control \cite{heemels2013a,heemels2013b,potoyan2013} and dynamic event-triggered control\cite{girard2014a} have been proposed to reduce the sensing effort to evaluate the event-triggered condition as well as to alleviate communication load for the NCSs. Moreover, model predictive control that incorporates the event/self-triggered control has been investigated in many works of literature, see, e.g.,\cite{hashimoto2018c,hashimoto2017c,hashimoto2019b,yulong,dai,dai2}. For more different formulations and approaches, see \cite{eventsurvey} for a recent survey paper. 

In this paper, we aim at investigating the event-triggered control strategy from different perspectives from the above previous works of literature. In particular, assuming that the dynamics of the plant is \textit{unknown} apriori, we consider the problem of \textit{finding} a set of suitable parameters for the event-triggered condition, such that certain specifications on the convergence of states and the safety are guaranteed. Finding such parameters are quite useful in practice to design the event-triggered controller, especially for the case when an accurate model of the plant is hard to obtain based on the first principles from physics, due to the fact that the dynamics is complex and highly nonlinear. 
To formulate our approach, we first define the indices that represent the convergence and safety properties that we want to evaluate. 
The goal is then to find a set of parameters for the event-triggered condition, such that some conditions on these indices are satisfied. Unfortunately, these indices are, when represented as the mappings from the parameters for the event-triggered condition, \textit{unknown} functions due to the fact that the dynamics of the plant is unknown apriori. Hence, we estimate these indices based on the training data, in particular employing the Gaussian process (GP) regression. 
The GP regression offers many benefits, such as the ability to incorporate prior knowledge about the model (e.g., smoothness, periodicity) by selecting suitable kernel functions, as well as the ability to provide uncertainty of the model for prediction values. In particular, in order to explore the desired parameter space, we employ the approach based on the \textit{safe} active learning, in which, for each iteration, the parameter with having the largest variance is evaluated. As we will see later, this approach is useful for {safety critical systems}, where violating the safety should not arise during the exploration. 
As for the theoretical analysis, we will mainly investigate whether the derived parameter space satisfies both the convergence and the safety specifications. 
In particular, the result shows that, under certain smoothness assumptions on the convergence and the safety indices, the exploration strategy succeeds to obtain the parameter space for the event-triggered controller, such that the resulting state trajectory satisfies both the convergence and the safety specifications. 

In summary, the contribution of this paper is: 
\begin{enumerate}
\item We propose a safe active learning framework to find a parameter space for the event-triggered condition. In particular, this parameter space is found, such that both the convergence and the safety specifications are satisfied; 
\item We provide a theoretical analysis to prove that the derived parameter space achieves both the convergence and the safety specifications. 
\end{enumerate}

So far, some model-free techniques to design the event-triggered control have been proposed in the literature, see, e.g., \cite{modelfree1,deepev,modelfree2}. For example, in \cite{modelfree1,modelfree2}, an actor-critic based $Q$-learning algorithm was proposed to learn the intermittent feedback controller. 
The approach presented in this paper differs from those previous works, in the sense that: (i) we provide an active learning framework based on the GP regression, aiming at finding the parameter space for the event-triggered condition; 
(ii) we can derive the parameter space that guarantees not only the convergence but also the safety specifications; 
(iii) it achieves a \textit{safe} exploration, meaning that the safety specification is always satisfied even in the exploration phase. 



\smallskip
\textit{Notation.} Let $\mathbb{N}$, $\mathbb{N}_{\geq 0}$, $\mathbb{N}_{>0}$ be the set of integers, non-negative integers, and positive integers, respectively. 
Let $\mathbb{R}$, $\mathbb{R}_{\geq 0}$, $\mathbb{R}_{>0}$ be the set of reals, non-negative reals, and positive reals, respectively. We denote by $\| \cdot \|$ the Euclidean norm. 
\section{Problem formulation} 
We consider the following nonlinear continuous-time systems: 
\begin{align}\label{dynamics}
\dot{x}(t) = f (x(t), u(t)), \ x(0) = \bar{x}, 
\end{align}
for all $t \in \mathbb{R}_{\geq 0}$, where $x \in \mathbb{R}^{n_x}$ is the state, $u \in \mathbb{R}^{n_u}$ is the control input, $\bar{x}$ is the initial state, and $f: \mathbb{R}^{n_x} \times \mathbb{R}^{n_u} \rightarrow \mathbb{R}^{n_x}$ is the map representing \textit{unknown} transition dynamics. 
Although the transition dynamics is unknown, it is assumed that the equilibrium point is known apriori;  without loss of generality, we assume that the equilibrium point is the origin, i.e., $0 = f(0,0)$. The control goal is to stabilize the state towards the origin. 

\begin{figure} [t]
	\centering
	\includegraphics[width = 5.0cm]{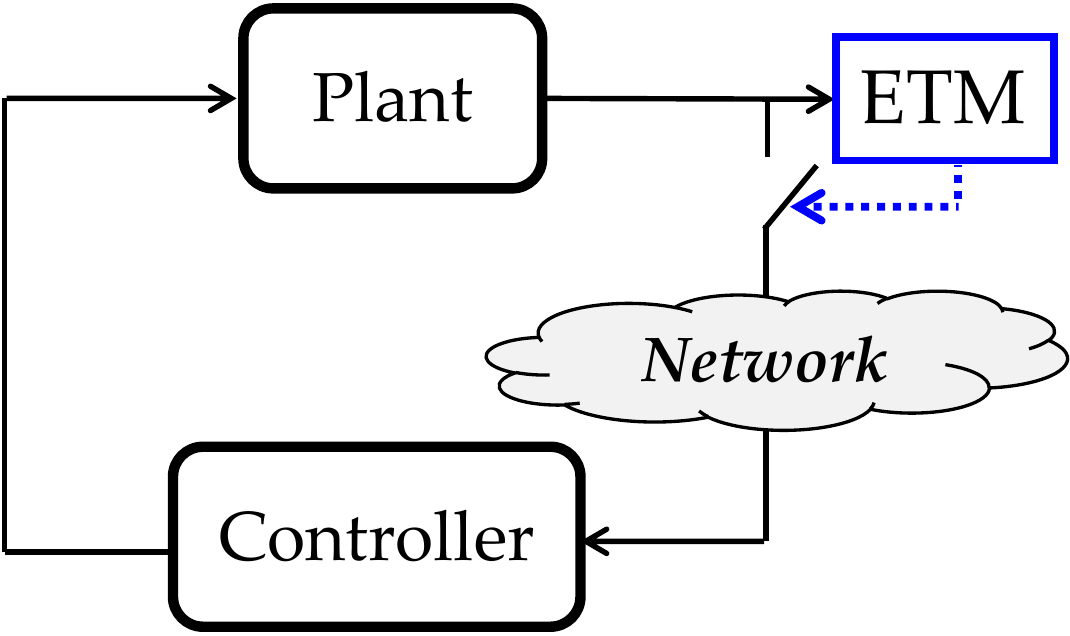}
	\caption{Event-triggered control system.}
	\label{NCS}
\end{figure}

In this paper, we aim to design an event-triggered control as shown in \rfig{NCS}. As shown in the figure, the event-triggered mechanism (ETM) is equipped along with the sensor, which determines the communication time to transmit the state measurement to the controller over the communication network. To illustrate the role of the ETM in more detail, let $t_k$, $k\in\{0, 1, 2, ...\}$ be the communication time instants when the ETM transmits the state to the controller. These time instances are recursively determined based on the state measurements as follows: 
\begin{align}\label{communicationtime}
t_{k+1} = {\inf} \{ t > t_k : h(x(t), x(t_k); \theta) \geq 0 \}, 
\end{align}
with $t_0 = 0$, where $h : \mathbb{R}^{n_x}\times \mathbb{R}^{n_x} \rightarrow \mathbb{R}_{\geq 0}$ is a function that characterizes the event-triggered condition and, as shown in \req{communicationtime}, it is parameterized by $\theta \in \mathbb{R}^{n_{\theta}}$ that we aim to design. 
Example of this function includes: 
\begin{align}
h(x(t), x(t_k); \theta) = \|x(t) - x(t_k) \| - \theta \|x(t)\| 
\end{align}
where $\theta \in \mathbb{R}_{\geq 0}$ (i.e., $n_\theta = 1$). 
Moreover, if the dynamic event-triggered controller\cite{girard2014a} is given, the parameter space will be multi-dimensional (i.e., $n_{\theta} > 1$). 
It is assumed that $\theta \in {\cal K} \subset \mathbb{R}^{n_\theta}$, where ${\cal K}$ is a given parameter space and is assumed to be compact. 
The controller is assumed to be given on the form: 
\begin{align}\label{controller}
u(t) & = \pi (\hat{x}(t)), 
\end{align}
where $\pi (\cdot)$ is the control law, and $\hat{x}(t)$ denotes the sampled state given by $\hat{x}(t) = x(t_k), \ \forall t \in [t_{k}, t_{k+1})$. For simplicity, in this paper it is assumed that the control law $\pi (\cdot)$ has been \textit{already} designed based on some model-free techniques under the time-triggered control. For example, we can design the PID controller such that the state is stabilized to the origin under the time-triggered controller. The controller tuning can be done by several techniques, such as the Bayesian optimization\cite{bo2}. 

Now, let ${\bf x} =\{x(t) \in \mathbb{R}^{n_x}, t\in\mathbb{R}_{\geq 0}\}$ denote the state trajectory following the dynamics of \req{dynamics}. 
Based on this trajectory, we proceed by defining the so-called \textit{convergence index} and the \textit{safety index}, which are denoted as $g({\bf x})$ and $s({\bf x})$, respectively. 
In particular, we say that the trajectory ${\bf x}$ satisfies the \textit{convergence specification} if $g({\bf x}) >0$. Moreover, we say that ${\bf x}$ satisfies the \textit{safety specification} if $s({\bf x}) >0$. For example, one can select the convergence index as 
\begin{align}\label{convergence}
g({\bf x}) = \inf_{t>0} \left \{ \cfrac{\eta (t)}{x(t) ^\mathsf{T} Q x(t)} -1 \right \}
\end{align}
where $Q$ is a given positive definite matrix, and $\eta(t)$ is characterized as $\eta(t) = \eta_0 \exp (-\beta t)$, so that $\eta(t) \rightarrow 0$, $t\rightarrow \infty$. 
Clearly, if $\inf_{t>0} ({\eta (t)}/{x(t) ^\mathsf{T} Q x(t)}) >1$, we have $\eta(t) > {x(t) ^\mathsf{T} Q x(t)}$, $\forall t\in \mathbb{R}$, so that $x(t) \rightarrow 0, t\rightarrow \infty$. Hence, if $g({\bf x}) > 0$, we have $x(t) \rightarrow 0, t\rightarrow \infty$ and so the convergence to the origin is guaranteed. Moreover, an example of the safety index includes
\begin{align}\label{safety}
s({\bf x}) = \inf_{t>0} \left \{ \xi - \|x(t)\| \right \}, 
\end{align}
where $\xi >0$ is a given threshold that characterizes the safety. 
That is, if $s({\bf x}) > 0$, we have $\|x(t)\| < \xi$, $\forall t\in\mathbb{R}$. As will be seen below, we aim at investigating the event-triggered parameter, such that both the convergence and the safety specifications are satisfied. 

Note that, since the controller is given, the state trajectory ${\bf x}$ is determined through the choice of the parameters for the event-triggered strategy, i.e., $\theta \in {\cal K}$. 
Hence, we can redefine the convergence and the safety indices as the mappings from $\theta$ as follows: 
\begin{align}
\widetilde{g} : {\cal K} \rightarrow \mathbb{R},\ \ \ \widetilde{s} : {\cal K} \rightarrow \mathbb{R}.  
\end{align}
That is, $\widetilde{g}$ and $\widetilde{s}$ represent the mappings from each $\theta \in {\cal K}$ onto the corresponding values of the convergence index and the safety index, respectively. 
Note that both $\widetilde{g}$ and $\widetilde{s}$ are the \textit{unknown} functions, due to the fact that the dynamics of the plant is unknown. 
Hence, as will be described in later sections, these functions are estimated by solving a certain regression problem based on the input-output training data: 
\begin{align}\label{trainingdata}
{\cal D} _g = \{\theta_i , y^g _i \}_{i\in\mathbb{N}_{>0}},\ \ 
{\cal D}_s = \{\theta_i , y^s _i \}_{i\in\mathbb{N}_{>0}}
\end{align} 
where $\theta_i $, $i\in\{1, 2, ...\}$ are the training inputs, and $y^g _i $ and $y^s _i$ are the corresponding (noisy) outputs of the convergence and the safety indices, respectively, and are given by 
\begin{align}
y^g _i = \widetilde{g}(\theta_i) + \epsilon^g _i,\  y^s _i = \widetilde{s}(\theta_i) + \epsilon^s _i. 
\end{align} 
In the above, $\epsilon^g _i$ and $\epsilon^s _i$ represent the additive noise, which are assumed to be uniformly chosen as $|\epsilon^g _i| \leq \eta_g$, $|\epsilon^s _i| \leq \eta_s$ for given $\eta_g >0$, $\eta_s>0$. For instance, the additive noise may arise due to the numerical error of computing the convergence/safety index. 

Based on the above notations and definitions, in this paper we consider the following problem: 
\begin{mypro}\label{problem}
\textit{Find the parameter space $\Theta \subset {\cal K}$ given by 
\begin{align}
\Theta = \left\{\theta \in {\cal K} : \widetilde{g}(\theta) > 0\ {\rm and}\ \widetilde{s}(\theta) >0 \right\}.
\end{align} 
Moreover, while collecting the training data as in \req{trainingdata}, it must satisfy the safety specification, i.e.,  $\widetilde{s} (\theta_i ) > 0$, for all $i \in\{1,2, ...\}$. } \qedwhite 
\end{mypro}

That is, we aim at exploring the parameter space in ${\cal K}$, such that both convergence index and the safety index are satisfied. Moreover, while collecting the training data to estimate $\widetilde{g}$ and $\widetilde{s}$, we must carefully select the training inputs $\theta_i $, $i\in\{1,2, ...\}$ such that the corresponding outputs always satisfy the safety specification $\widetilde{s} (\theta_i ) > 0$, for all $i \in\{1,2, ...\}$. Such safety constraint should be crucially needed if we would like to avoid the violation of the safety during the experiment. 

\section{Safe and active exploration}
In this section we provide a solution approach to \rpro{problem}. 

\subsection{Main algorithm}
Recall that, since $\widetilde{g}$ and $\widetilde{s}$ are the unknown functions, we estimate them based on the training data as shown in \req{trainingdata}. 
In particular, in this paper we estimate these functions by the Gaussian process (GP) regression \cite{rasmussen}. 
Similarly to \req{trainingdata}, let ${\cal D}_{g, N} = \{\theta_i , y^g _i \}^ N_{i=1}$, ${\cal D}_{s, N} = \{\theta_i , y^s _i \}^N _{i=1}$ with $N$ being the number of the training data. 
Moreover, let $\mathsf{k}_g : {\cal K} \times {\cal K} \rightarrow \mathbb{R}_{\geq 0}$ be a given kernel function and ${K}_{g} \in \mathbb{N}^{N\times N}$ be the covariance matrix for the estimate of $\widetilde{g}$, such that $K^{(a, b)}_g = \mathsf{k}_g (\theta^{(a)} _{ev}, \theta^{(b)} _{ev})$, where $K^{(a, b)} _g$ denotes the $(a, b)$-component of $K_g$. 
Then, the GP regression estimates $\widetilde{g}$ as 
\begin{align*}
\widetilde{g}_N (\theta) \sim \mathcal{N} \left (\mu_{g,N}  \left(\theta\right), {\sigma}^2_{g, N} \left({{\theta}}\right)\right), 
\end{align*}
where $\widetilde{g}_N$ denotes the GP model of $\widetilde{g}$ based on the training data $\mathcal{D}_{g,N}$, and 
\begin{align}
\mu_{g, N}  \left(\theta\right) &= \mathsf{{k}}^\mathsf{T} _{*g} ({\theta}) ({K}_{g} + \eta^2 _g I_N) ^{-1} {Y} _{g, N} \notag \\
{\sigma}^2_{g, N} \left(\theta\right) &= \mathsf{k}_g (\theta , \theta) - \mathsf{k}^\mathsf{T} _{*g} ({ \theta}) ({K}_{g} + \eta^2 _g I_N) ^{-1} \mathsf{k}_{*g} (\theta), \notag
\end{align}
with ${Y}_{g, N}, \mathsf{k} _{*g}(\theta) \in \mathbb{R}^{N}$ being given by
\begin{align}
&{Y}_{g, N} = [y^g _1,..., y^g _N]^\mathsf{T}, \\
&\mathsf{k} _{*g} ({ \theta}) = \left [\mathsf{k}_g ({\theta} ,\theta_1), ..., \mathsf{k}_g (\theta , {\theta_N}))\right]^\mathsf{T}. 
\end{align}
Similarly, letting $\mathsf{k}_s : {\cal K} \times {\cal K} \rightarrow \mathbb{R}_{\geq 0}$ be a given kernel function and ${K}_{s} \in \mathbb{N}^{N\times N}$ be the covariance matrix for the estimate of $\widetilde{s}$, we have $\widetilde{s}_N (\theta) \sim \mathcal{N} \left (\mu_{s,N}  \left(\theta\right), {\sigma}^2_{s, N} \left({{\theta}}\right)\right)$
where $\mu_{s, N}  \left(\theta\right) = \mathsf{{k}}^\mathsf{T} _{*s} ({\theta}) ({K}_{s} + \eta^2 _s I_N)^{-1} {Y} _{s, N}$, ${\sigma}^2_{s, N} \left(\theta\right) = \mathsf{k}_s (\theta , \theta) - \mathsf{k}^\mathsf{T} _{*s} ({ \theta})({K}_{s} + \eta^2 _s I_N) ^{-1} \mathsf{k}_{*s} (\theta)$, 
with ${Y}_{s, N}$ and $\mathsf{k} _{*s} ({ \theta})$ being defined in the same way as for $\widetilde{g}$. Before providing the algorithm, we need to make the following assumption:  
\begin{myas}\label{as1}
\textit{It is known that there exist $\Theta_{init} \subset {\cal K}$, such that $\widetilde{s} (\theta) >0$, for all $\theta \in \Theta_{init}$.} \qedwhite 
\end{myas}

\ras{as1} implies that, we have the prior knowledge (before the exploration) that $\Theta_{init}$ is the parameter space guaranteeing the safety specification. 
This assumption is required to explore the parameter space at the initial phase and collect the training data while guaranteeing the safety specification. 

\begin{algorithm}[htbp]
{\small
\SetKwInOut{Input}{Input}
\SetKwInOut{Output}{Output}
\Input{Characterizations for $\beta_{s,N}$, $\beta_{g,N}$, $N=1, \ldots$ (see \rsec{thsec}); $\Theta_{init}$ (initial parameter space); } 
\Output{$\Theta$ (parameter space satisfying convergence and safety);}
${\cal D}_g \leftarrow \varnothing$, ${\cal D}_s \leftarrow \varnothing$; \\
\textit{[Initial exploration phase]:}\\ 
\For{$i=1:N_{init}$}{
Randomly select $\theta \in \Theta_{\rm init}$. \\
Given $\theta = [\sigma_{0}, \sigma_{\infty}]$, implement the event-triggered strategy and measure the outputs $y_g$, $y_s$. \\
Set the training data as $\mathcal{D}_g \leftarrow \mathcal{D}_g \cup \{\theta, y_g\}$, $\mathcal{D}_s \leftarrow \mathcal{D}_s \cup \{\theta, y_s\}$; 
}
$N \leftarrow N_{init}$;  \\
Estimate $\widetilde{g}$ and $\widetilde{s}$ based on the GP regression: 
\begin{align}
\widetilde{g}_N (\theta) &\sim \mathcal{N} \left (\mu_{g,N}  \left(\theta\right), {\sigma}^2_{g, N} \left({{\theta}}\right)\right)\\ 
\widetilde{s}_N (\theta) &\sim \mathcal{N} \left (\mu_{s,N}  \left(\theta\right), {\sigma}^2_{s, N} \left({{\theta}}\right)\right)
\end{align}
for all $\theta \in {\cal K}$. \\

Set $\Theta, \Theta_s \leftarrow \Theta_{init}$; 

\smallskip
\textit{[Exploration phase]:} \\
\For{$j = 1:N_{exp}$} {
\textit{(Step~1)} Compute ${\Theta}_{s, N}, {\Theta}_N \subseteq {\cal K}$ as  
\begin{align}
{\Theta}_{s, N} = \{ \theta \in {\cal K} : & \mu_{s, N}(\theta) - \beta_{s, N} \sigma _{s, N} (\theta) > 0 \},\notag 
\end{align}
\begin{align}
{\Theta}_{N} = \{ \theta \in {\cal K} : &\mu_{s, N}(\theta) - \beta_{s, N} \sigma _{s, N} (\theta) > 0\ {\rm and} \notag \\
                 &\ \mu_{g, N}(\theta) - \beta_{g, N} \sigma _{g, N} (\theta) > 0. \notag 
\end{align}
Then, set $\Theta_s \leftarrow \Theta_s \cup {\Theta}_{s, N}$, $\Theta \leftarrow \Theta \cup {\Theta}_N$; 

\smallskip
\textit{(Step~2)} Select $\theta_i  \in {\cal K}$ by
\begin{align}
\theta_i  &=\underset{\theta\in\Theta_s }{\arg\max}\ \{\sigma^2 _{g, N} (\theta) +  \sigma^2 _{s, N} (\theta) \} \label{max1}. 
\end{align}

\smallskip
\textit{(Step~3)} Given $\theta_i $, implement the event-triggered strategy and measure the corresponding outputs $y^g _i$, $y^s _i$. Then, update the training data as 
\begin{align}
\mathcal{D}_g &\leftarrow \mathcal{D}_g \cup \{\theta_i , y^g _i \}, \\
\mathcal{D}_s &\leftarrow \mathcal{D}_s \cup \{\theta_i , y^s _i \}; 
\end{align} 

\textit{(Step~4)} Set $N \leftarrow N+1$ and update the estimate of $\widetilde{g}$ and $\widetilde{s}$ based on the GP regression: 
\begin{align}
\widetilde{g}_N (\theta) &\sim \mathcal{N} \left (\mu_{g,N}  \left(\theta\right), {\sigma}^2_{g, N} \left({{\theta}}\right)\right)\\ 
\widetilde{s}_N (\theta) &\sim \mathcal{N} \left (\mu_{s,N}  \left(\theta\right), {\sigma}^2_{s, N} \left({{\theta}}\right)\right)
\end{align}
}
Return $\Theta$. 
   \caption{Safe and active parameter exploration.}\label{exploration}
    }
\end{algorithm}

The approach presented in this paper is based on the GP-based active exploration (see, e.g., \cite{safeexploration2015a,sui2015a}), and the details are shown in \ralg{exploration}. In the algorithm, we let $N_{init}$ and $N_{exp}$ denote the number of iterations for the initialization and for the implementation of the exploration, respectively. The algorithm starts by exploring the initial parameter space that guarantees safety $\Theta_{init}$ (initial exploration phase), which is then followed by the main algorithm for the active exploration (exploration phase). 
As shown in the algorithm, the exploration phase implements the four steps (Step~1, 2, 3, 4). 
First, we compute the two sets $\Theta_s$, $\Theta$ (Step~1). 
Intuitively, $\Theta_s$ represent the the parameter space that guarantee the safety specification, and $\Theta$ represents the parameter space that guarantees both the safety and the convergence specifications. In the algorithm, $\beta_{g, N}, \beta_{s, N} >0$ are the weight parameters that specify the confidence bounds. A detailed characterization on how to select these parameters will be given in the next subsection, where we provide a theoretical analysis for proving the convergence and the safety specification. 
Then, we determine the event-triggered parameter $\theta_i$ by solving \req{max1} (Step~2). Intuitively, we pick up the parameter from $\Theta_s$ (i.e., the parameter space that guarantees safety) that provides the largest uncertainties about the unknown functions $\widetilde{g}$, $\widetilde{s}$ in order to effectively explore the parameter space while guaranteeing the safety specification. 
Once this parameter is chosen, we implement the event-triggered strategy to measure the corresponding output and update the training data (Step~3), and, finally, update the GP model of $\widetilde{g}$ and $\widetilde{s}$ (Step~4). 
If the above steps are iterated for $N_{exp}$ times, we return $\Theta$ as the parameter space that guarantees safety and convergence specifications. 

\subsection{Theoretical analysis}\label{thsec}
In this section we investigate a theoretical analysis to show that the derived parameter space $\Theta$ guarantees both the safety and the convergence specifications, and, at the same time it ensures the safe exploration, i.e., $\widetilde{s} (\theta_i) >0$, $\forall i\in\{1, 2, ...\}$. Before providing the analysis, we need to make certain \textit{smoothness} assumptions on the unknown functions $\widetilde{s}$, $\widetilde{g}$. In particular, it is assumed that both $\widetilde{g}$ and $\widetilde{s}$ lie in the reproducing kernel Hilbert space (RKHS) corresponding to the kernel $\mathsf{k}_g$ and $\mathsf{k}_s$, respectively. Moreover, it is assumed that $\| \widetilde{g} \|_{\mathsf{k}_g} \leq B_g$, $\| \widetilde{s} \|_{\mathsf{k}_s} \leq B_s$ for some $B_s, B_g>0$, where $\| \cdot \|_{\mathsf{k}_g}$ and $\| \cdot \|_{\mathsf{k}_s}$ denote the induced norms of the RKHS. The assumption that $\widetilde{g}$ has a bounded norm in the RKHS implies that the function is given of the form $\widetilde{g} (\theta) = \sum^\infty _{n=1} \alpha_n \mathsf{k}_g (\theta, \theta^* _{n})$, where $\theta^* _{n} \in {\cal K}$, $n\in\mathbb{N}_{>0}$ are the representer points and $\alpha_n \in \mathbb{R}$, $n\in\mathbb{N}_{>0}$ are the parameters that decay sufficiently fast as $n$ increases, so that $\sum^\infty _{n=1} \alpha_n < \infty$ ($\widetilde{s} (\theta)$ can be defined in a similar manner). 
The following result shows that, under the above assumptions as well as the suitable selections on $\beta_{s, N}, \beta_{g, N}$, the derived parameter space is guaranteed to satisfy the convergence and the safety specifications, as well as satisfy the safe exploration. 
\begin{mythm}\label{result}
\textit{Suppose that $\widetilde{g}$ and $\widetilde{s}$ lie in the reproducing kernel Hilbert space (RKHS) corresponding to the kernel $\mathsf{k}_g$ and $\mathsf{k}_s$, respectively, and that $\| \widetilde{g} \|_{\mathsf{k}_g} \leq B_g$, $\| \widetilde{s} \|_{\mathsf{k}_s} \leq B_s$. 
Moreover, let $\beta_{g, N} > 0$, $\beta_{s, N} > 0$ be given by 
\begin{align}
\beta_{g, N} &=  \sqrt{B^2 _g - {Y}^\mathsf{T} _{g, N} K_{g}^{-1}{Y} _{g, N} + N } \\
\beta_{s, N} &=  \sqrt{B^2 _s - {Y}^\mathsf{T} _{s, N} K_{s}^{-1}{Y} _{s, N} + N } 
\end{align}
for all $N = 1, 2, ...$. Then, running Algorithm~1 ensures that: 
\begin{enumerate}
\item Let $\Theta \subseteq {\cal K}$ be the resulting output from Algorithm~1. 
Then, $\widetilde{s}(\theta) > 0$, $\widetilde{g}(\theta) >0$ for all $\theta \in \Theta$; 
\item While implementing Algorithm~1, it ensures the safe exploration, i.e., $\widetilde{s} (\theta_i) >0$, $\forall i\in\{1, 2, ...\}$.  \qedwhite
\end{enumerate}}
\end{mythm} 

The proof is based on \cite{srinivas} and the overview is described as follows. 

\begin{proof}
\textit{(Sketch)} We have 
\begin{align}
|{\mu}_{g, N} (\theta) - \widetilde{g} (\theta)|  &\leq \mathsf{k}_{g, N}(\theta, \theta)^{-1/2} \|{\mu}_{g} - \widetilde{g}\|_{\mathsf{k}_{g, N}} \label{first} \\ 
&= {\sigma}_{g, N} (\theta) \|{\mu}_{g, N} - \widetilde{g}\|_{\mathsf{k}_{g, N}}, \notag 
\end{align}
where the first inequality follows from the Cauchy-Schwarz inequality, and $\mathsf{k}_{g, N} : {\cal K} \times {\cal K} \rightarrow \mathbb{R}_{\geq 0}$ is given by 
\begin{align*}
\mathsf{k}_{g, N}(\theta , \theta') = &\mathsf{k}_g (\theta , \theta' ) - \mathsf{k}^\mathsf{T} _{*g} (\theta) {K}_{g}^{-1} \mathsf{k}_{*g} (\theta'). 
\end{align*}
Moreover, in \req{first}, $\|\cdot \|_{\mathsf{k}_{s, N}}$ denotes the norm induced by the RKHS corresponding to $\mathsf{k}_{s, N}$. 
Then, from the proof of Lemma~7.2 in \cite{srinivas}, it follows that
\begin{align}
\|{\mu}_{g, N} - \widetilde{g} \|^2 _{\mathsf{k}_{g, N}} & = B^2 _g - {Y}^\mathsf{T} _{g, N} K_{g} ^{-1}{Y} _{g, N} + \eta^{-2} _g \sum^N _{i=1} (\epsilon^g _{i})^2 \notag \\ 
& \leq B^2 _g - {Y}^\mathsf{T} _{g, N} K_{g} ^{-1}{Y} _{g, N} + N, 
\end{align}
where we have used $|\epsilon^g _{i}| \leq \eta_g$, $i\in 1, 2, ...$. 
Hence, letting $\beta_{N, g} = \sqrt{B^2 _g - {Y}^\mathsf{T} _{g, N} K_{g} ^{-1}{Y} _{g, N} + N}$, we have 
\begin{align}
|{\mu}_{g, N}(\theta) - \widetilde{g} (\theta)| \leq \beta_N {\sigma}_{g, N} (\theta), 
\end{align}
for all $\theta \in {\cal K}$. 
Therefore, for all $\theta \in {\cal K}$, 
\begin{align}
{\mu}_{g, N}(\theta) &- \beta_N {\sigma}_{g, N} (\theta) \notag \\ 
&= {\mu}_{g, N}(\theta) - \widetilde{g}(\theta) - {\sigma}_{g, N} (\theta) + \widetilde{g}(\theta), \notag \\ 
& \leq \widetilde{g}(\theta)
\end{align}
and so ${\mu}_{g, N}(\theta) - \beta_{g, N} {\sigma}_{g, N} (\theta) > 0$ implies $\widetilde{g}(\theta) > 0$. 
Similarly, ${\mu}_{s, N}(\theta) - \beta_{s, N} {\sigma}_{s, N} (\theta) > 0$ implies $\widetilde{s}(\theta) > 0$. 
Hence, letting $\Theta = \bigcup^{N_{ini} + N_{exp}} _{N=1} \Theta_N\subseteq {\cal K}$, where $\Theta_N$ is given by 
\begin{align}
\Theta_N = \{\theta \in {\cal K}: &\ \mu_{g, N} (\theta) - \beta_{g, N} (\theta) >0\ {\rm and} \notag \\ 
 &\ \ \mu_{s, N} (\theta) - \beta_{s, N} (\theta) >0 \}, 
\end{align} 
it then follows that $\widetilde{s} (\theta) > 0$, $\widetilde{g} (\theta) > 0$, $\forall \theta \in \Theta$. 
Hence, the derived parameter space $\Theta$ satisfies the convergence and the safety specifications. 
The safe exploration property $\widetilde{s}(\theta_N) > 0$ for all $N = 1, 2, ..., N_{exp}$ can be similarly proven, since, for each iteration $j$, we select $\theta_j$ from the parameter space that guarantees the safety specification $\Theta_s$ (see \req{max1}). 
\end{proof}

\section{Simulation result} 
Consider the following inverted pendulum system: 
\begin{align}
\dot{x}_1(t) &= x_2(t) \\
\dot{x}_2(t) &= \sin x_1(t) - x_2(t) + u(t),  
\end{align}
with $x_0 = [1.0,\ 0]^\mathsf{T}$ and $u \in \mathbb{R}$ is the control input. Note that the above dynamics is assumed to be \textit{unknown} apriori. The control law has been designed by $u(t) = K \hat{x}(t)$ with $K= [-1.08,\ -1.43]$, which has been found based on the Bayesian optimization technique, see, e.g., \cite{bo2}; the detailed procedure is omitted for brevity. It can be verified that, under the time-triggered controller (i.e., the control input is updated continuously) with the derived control gain $K$, the state is stabilized to the origin, i.e., $x(t) \rightarrow 0$, $t\rightarrow \infty$. The function representing the event-triggered strategy is given by 
\begin{align}
h(x(t), x(t_k); \theta) = \|x(t) - x(t_k) \| - \varepsilon (t) \|x(t)\| 
\end{align}
where $\varepsilon (t)$ is a time-varying parameter characterized as 
\begin{align}\label{epsilon}
\varepsilon (t) = (\varepsilon_{0} - \varepsilon_{\infty} ) \exp(-\gamma t) + \varepsilon_{\infty} 
\end{align}
for all $t\in\mathbb{R}$, where $\gamma \in \mathbb{R}_{> 0}$ is the parameter representing the rate for the convergence of $\varepsilon(t)$, and $\varepsilon_{0}, \varepsilon_{\infty} \in \mathbb{R}_{\geq 0}$ represent thresholds for the initial time and for the steady state time, respectively. The event-triggered strategy based on the above time-varying threshold is related to the dynamic event-triggered control (see, e.g., \cite{girard2014a}). 
For simplicity, it is assumed that $\gamma = 0.1$ is a \textit{given} parameter and $\varepsilon_{0}, \varepsilon_{\infty}$ are the parameters to be designed. In other words, we have $\theta = [\varepsilon_{0}, \varepsilon_{\infty}] \in \mathbb{R}^2 _{\geq 0}$. Moreover, the parameter space is given by ${\cal K} = [0.01, 1] \times [0.01, 1]$. 
The convergence index is given by \req{convergence} with $Q = I_2$ and $\eta(t) = \eta_0 \exp (-\beta t)$, $\eta_0 = 2.0$, $\beta=0.05$. 
The safety index is given by $s({\bf x}) = \inf_{t>0} \left \{ \xi - |x_2(t)| \right \}$ with $\xi = 0.25$ (i.e., the amount of the velocity should be less than $0.25$). The initial parameter space is assumed to be given by $\Theta_{init} = [0.01, 0.05] \times [0.01, 0.05]$. The number of iteration for Algorithm~1 is $N_{init} = 10$, $N_{exp} = 100$. 

\begin{figure} [t]
	\centering
	\includegraphics[width = 7.0cm]{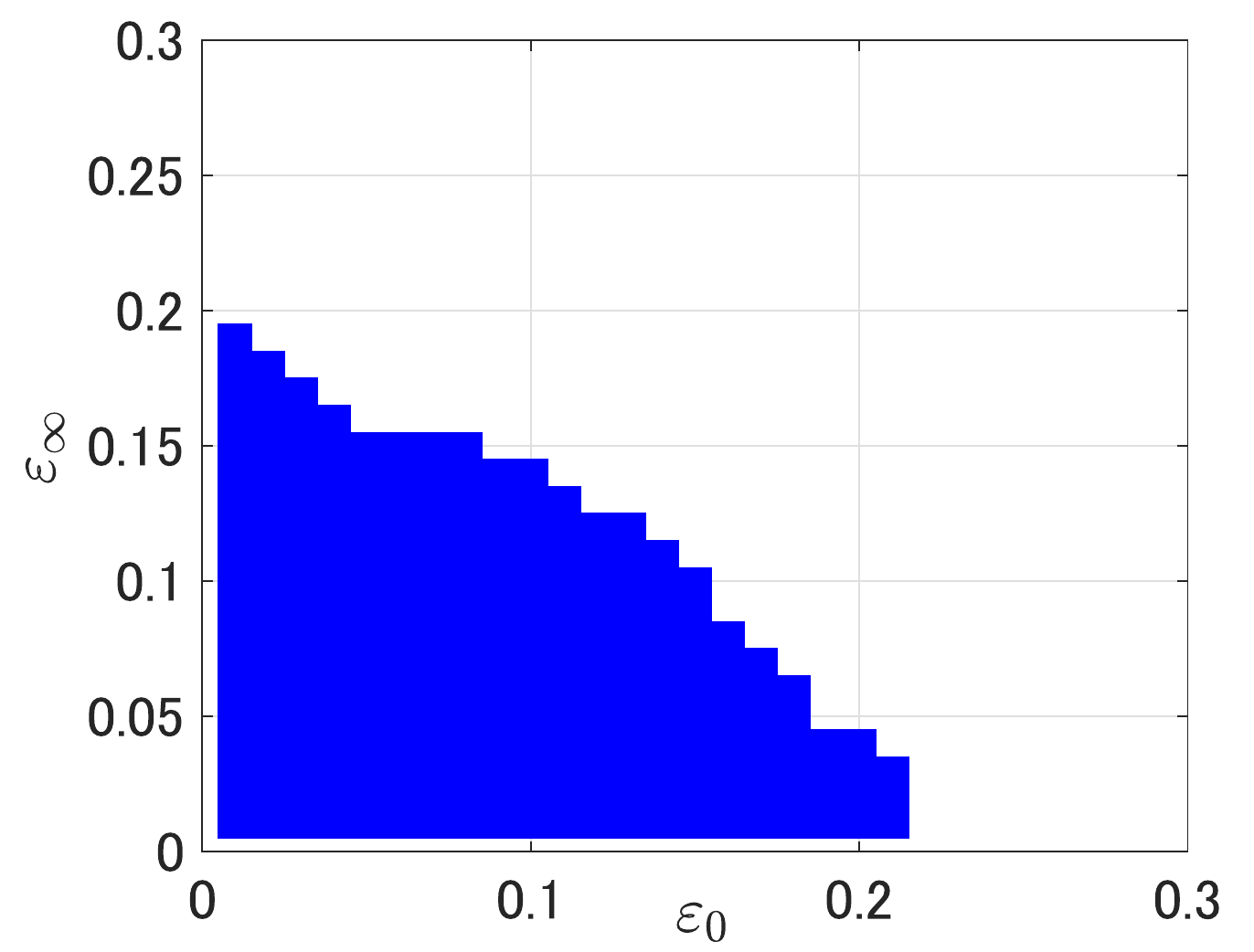}
	\caption{The parameter space $\Theta$ derived by applying Algorithm~1 (the blue region). }
	\label{param}
\end{figure}

\begin{figure} [t]
	\centering
	\includegraphics[width = 8.5cm]{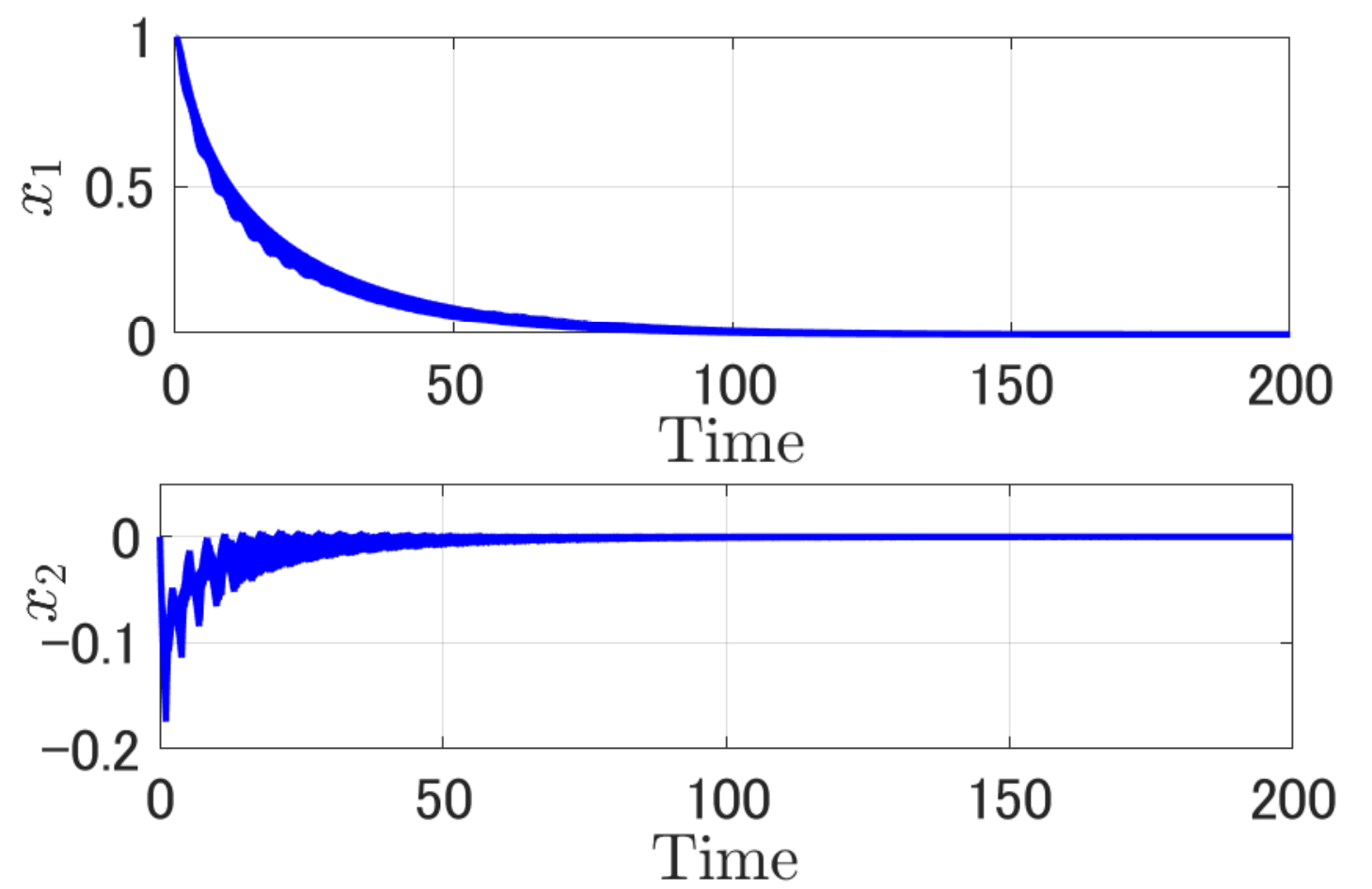}
	\caption{State trajectories with the parameter $\theta$ randomly chosen from $\Theta$. It can be verified that all state trajectories satisfy the convergence and the safety specifications.}
	\label{random}
\end{figure}

\begin{figure} [t]
	\centering
	\includegraphics[width = 8.0cm]{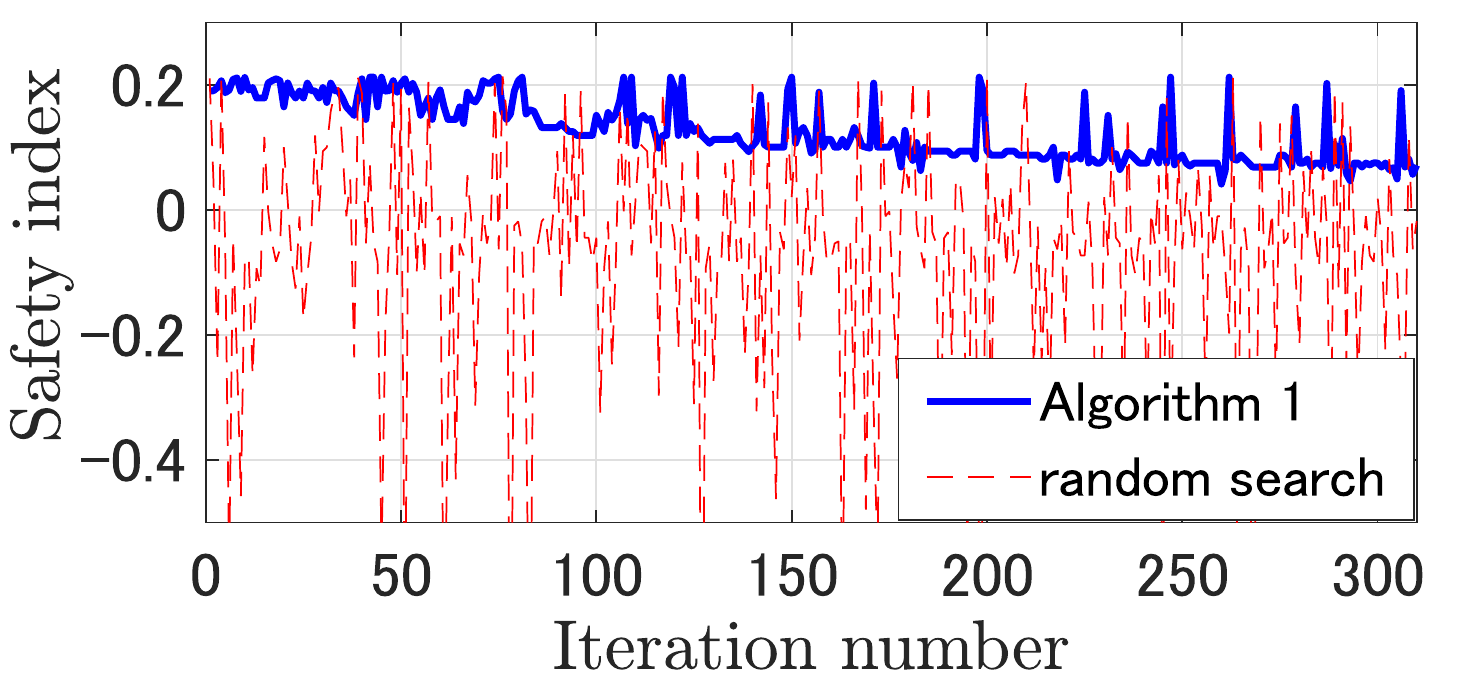}
	\caption{Safety index computed for each iteration number $j$ under Algorithm~1 (blue solid) and the random search (red dotted).}
	\label{safeindex}
\end{figure}

\rfig{param} represents the parameter space $\Theta$ obtained by applying Algorithm~1. 
Moreover, \rfig{random} illustrates the $100$ state trajectories with the event-triggered parameter randomly chosen from the computed $\Theta$. 
It can be verified that all state trajectories satisfy both convergence and the safety specifications, showing the effectiveness of the proposed approach. In particular, it can be seen that the amount of the velocity (i.e., $|x_2|$) is always less than $0.25$ and thus the safety specification is satisfied. 
In addition, \rfig{safeindex} illustrates the values of the safety index computed for each iteration number $j$ during the exploration phase. 
To make a comparison, we also illustrate the result of the \textit{random search}, where the parameter is uniformly chosen from the interval $[0.01, 1.0] \times [0.01, 1.0]$ (instead of selecting the parameter from \req{max1}). The figure shows that the random search fails to achieve the safe exploration, since the safety indices become negative for some $j$. On the other hand, our approach succeeds to achieve the safe exploration, since the safety indices are positive for all times during the implementation of Algorithm~1. 

\section{Conclusion and future work}
In this paper, we investigate a parameter exploration strategy for the event-triggered control, such that both the prescribed convergence and the safety specifications are satisfied. 
Moreover, we provide a theoretical analysis, such that, under the smoothness assumption on the convergence and the safety indices, the derived parameter space achieves both convergence and safety specifications. Future work involves investigating the applicability of our approach to practical implementations.

\end{document}